\documentclass[11pt,a4paper]{article}

\usepackage{inputenc}
\usepackage{amsmath}
\usepackage{bm}
\usepackage{bbold}
\usepackage{amsthm}
\usepackage{enumerate}

\usepackage[top=1.0in, left=1.0in, right=1.0in, bottom=1.0in]{geometry}

\usepackage{tikz}\tikzset{x=1cm,y=1cm,z=1cm}

\usepackage{pgfplots}\pgfplotsset{compat=1.16}

\usepackage[hyphens]{url}

\usepackage{hyperref}
\usepackage{breakurl}

\title{Algebraic solution of tropical best approximation problems\thanks{Mathematics 2023, 11(18), 3949; https://doi.org/10.3390/math11183949}}

\author{N. Krivulin\thanks{Faculty of Mathematics and Mechanics, Saint Petersburg State University, 28 Universitetsky Ave., St.~Petersburg, 198504, Russia, 
nkk@math.spbu.ru.}}

\date{}

\newtheorem{theorem}{Theorem}
\newtheorem{lemma}[theorem]{Lemma}

\theoremstyle{definition}
\newtheorem{algorithm}{Algorithm}

\begin{document}

\maketitle

\begin{abstract}
We introduce new discrete best approximation problems, formulated and solved in the framework of tropical algebra, which deals with semirings and semifields with idempotent addition. Given a set of samples, each consisting of the input and output of an unknown function defined on an idempotent semifield, the problem is to find a best approximation of the function, by tropical Puiseux polynomial and rational functions. A new solution approach is proposed, which involves the reduction of the problem of polynomial approximation to the best approximate solution of a tropical linear vector equation with an unknown vector on one side (a one-sided equation). We derive a best approximate solution to the one-sided equation, and we evaluate the inherent approximation error in a direct analytical form. Furthermore, we reduce the rational approximation problem to the best approximate solution of an equation with unknown vectors on both sides (a two-sided equation). A best approximate solution to the two-sided equation is obtained in numerical form, by using an iterative alternating algorithm. To illustrate the new technique developed, we solve example approximation problems in terms of a real semifield, where addition is defined as maximum and multiplication as arithmetic addition (max-plus algebra), which corresponds to the best Chebyshev approximation by piecewise linear functions.
\\

\textbf{Key-Words:} tropical semifield, tropical Puiseux polynomial, best approximate solution, discrete best approximation, Chebyshev approximation.
\\

\textbf{MSC (2020):} 15A80, 90C24, 41A50, 41A65, 65D15
\end{abstract}

\section{Introduction}

The discrete approximation problem --- where an unknown real function $y=f(x)$ is approximated, given a set of samples $(x_{i},y_{i})$ of function values $y_{i}=f(x_{i})$ at some points $x_{i}$ --- is of primary interest in data analysis. Let $F(x;\bm{\theta})$ be an approximating function that depends on the vector $\bm{\theta}$ of the unknown parameters to be determined. Consider a best approximate solution to the problem, which is defined in the sense of a distance function, $d$, used for estimating the approximation error, to find
\begin{equation}
\bm{\theta}_{\ast}
=
\arg\min_{\bm{\theta}}
\max_{i}d(F(x_{i};\bm{\theta}),y_{i}).
\label{P-thetaast}
\end{equation}

The minimax optimization problem at \eqref{P-thetaast} is commonly formulated in terms of the distance between two reals $a$ and $b$, given by the absolute value $d(a,b)=|a-b|$, which leads to the Chebyshev metric to measure the error. The problem is then referred to as the best Chebyshev approximation problem \cite{Harter1975Method,Mhaskar2000Fundamentals,Steffens2006History}, and it has drawn much attention over the centuries since Laplace's classical work \cite{Laplace1832Mecanique} (Book~3, Chapter~V, \S39). 

In many Chebyshev approximation problems, the approximating function $F(x;\bm{\theta})$ is taken as linear, with respect to parameters. These problems appear in the solution of overdetermined systems of linear equations \cite{Tewarson1972Minimax,Appa1973L1,Pinar2009Overdetermined}, regression analysis when the distribution of errors has bounded support \cite{Rabinowitz1968Applications,Appa1973L1,Hand1978Aspects,Kennedy1980Statistical} and in other applications in science and engineering.

Under the above linearity assumption, the problems can be reformulated as linear programs and then numerically solved by computational procedures and techniques available in linear programming. Examples of linear programming and other algorithmic solutions to the best Chebyshev approximation problem can be found in early works \cite{Wagner1959Linear,Stiefel1960Note,Rabinowitz1968Applications,Tewarson1972Minimax,Appa1973L1,Watson1973Best,Sposito1976Minimizing,Armstrong1979Algorithm} and in more recent publications \cite{Kim2000Algorithms,Boyd2004Convex,Castillo2008Dealing,Ene2019Improved,Krivulin2020Usingparameter}.

If the approximating function $F(x;\bm{\theta})$ is assumed to be nonlinear in its parameters, the solution of the best Chebyshev approximation problem becomes more difficult, since it involves the minimization of a complicated nonconvex and nondifferentiable function. In this case, the problems are handled by appropriate methods and techniques of nonlinear programming \cite{Conn1988Computational}, which may use various gradient and gradient-free optimization schemes. As typical solutions, one can consider iterative computational procedures that solve an auxiliary linear \cite{Osborne1969Algorithm,Madsen1975Minimax,Anderson1977Discrete,Charalambous1978Efficient} or quadratic \cite{Hald1981Combined,Han1981Variable,Fletcher1982Model} programming problem at each iteration.

In this paper, we offer a new treatment of the best approximation problem, where the problem is put in the framework of tropical algebra, to develop new approximation models and solutions. Tropical algebra deals with the theory and methods of semirings and semifields with idempotent addition \cite{Kolokoltsov1997Idempotent,Golan2003Semirings,Heidergott2006Maxplus,Itenberg2007Tropical,Gondran2008Graphs,Butkovic2010Maxlinear,Maclagan2015Introduction} and is being increasingly applied in a range of fields, from algebraic geometry to decision making. An example of the tropical semifield, which is usually referred to as max-plus algebra, is an extended set of reals, where the addition is defined as maximum and the multiplication as arithmetic addition.

Application of tropical algebra allows one to represent various nonlinear problems as linear ones in the tropical algebra setting, which simplifies the solution of the problems and facilitates the interpretation of results. In many cases, the solution of problems in the framework of tropical algebra can lead to direct analytical results that describe all solutions in an explicit compact vector form.

Models and methods of tropical algebra are used in approximation problems in various settings, including problems of best approximation, with respect to the Hilbert projective metric \cite{Akian2011Best,Saadi2021Zerosum} and the Chebyshev metric \cite{Butkovic2010Maxlinear}. The known results mainly concern best approximation of a vector by elements of a tropical linear space (tropical linear regression), defined in the framework of max-plus algebra. The discrete best approximation problems, where given data are to be fitted with tropical polynomial and rational functions, present another promising direction of tropical approximation. These problems have appeared in a few recent works on neural networks and machine learning (see, e.g., \cite{Zhang2018Tropical,Maragos2021Tropical}) and require further development of related theory and solution techniques.

We consider approximation problem \eqref{P-thetaast}, reformulated to approximate functions defined on a tropical semifield (a semiring with idempotent addition and invertible multiplication). The approximation error is defined in the form of a generalized metric on the tropical vector space over the semifield. As the approximating function, we take tropical analogues of Puiseux polynomial and rational functions. We observe that in the case of max-plus algebra, the generalized metric coincides with the Chebyshev metric, whereas the Puiseux polynomial and rational functions are real convex and nonconvex piecewise linear functions. This leads to a natural application of the problem under consideration to the best Chebyshev piecewise linear approximation of real functions.

In conventional algebra, the Puiseux polynomials extend the notion of standard polynomials, by allowing rational exponents. Tropical Puiseux polynomials are introduced in the same way, where addition, multiplication and exponentiation are defined in terms of an idempotent semifield. These polynomials arise as both an effective instrument and the key object of analysis in a range of research contexts, from optimization problems in operations research to tropical algebraic geometry \cite{Itenberg2007Tropical,Markwig2010Field,Maclagan2015Introduction,Hampe2017Tropical,Grigoriev2018Tropical,Krivulin2021Algebraic}. Thus, the development of approximation techniques using tropical Puiseux polynomial and rational functions can be considered to be of benefit to both tropical algebra and its applications. 
 
The key ingredient of the proposed approximation technique is to turn the best discrete approximation problem into that of solving tropical linear vector equations with an unknown vector on one side of the equation (one-sided equations) or unknown vectors on both sides (two-sided equations). We reduce the approximation by tropical polynomials to finding a best approximate solution of a one-sided vector equation, and the approximation by tropical rational functions to a best approximate solution of a two-sided equation. We handle the one-sided equation by applying the results in \cite{Krivulin2009Solution,Krivulin2013Solution-linear}, which offer a direct explicit analytical solution to the problem. To obtain a best approximate solution of the two-sided equation, we use the iterative alternating algorithm proposed in \cite{Krivulin2023Solution}.

As the main contribution of the study, we consider the development of a new analytical framework for the best approximation of functions in the tropical algebra setting. We propose a new formulation of the approximation problem, find a strong connection between the approximation problem and tropical linear vector equations and, finally, provide a workable solution to the problem. We demonstrate the applicability of the proposed technique, with examples of best approximation problems formulated in the max-plus algebra setting. Both the problems of approximating a convex function by tropical polynomials and a nonconvex function by tropical rational functions are considered. We obtain the best approximate solutions to the problems and give graphical illustrations.

Among the advantages of the proposed approach over the other solutions to the discrete best approximation problems in the tropical algebra setting is the formulation and solution of the problem in general terms of an arbitrary tropical semifield. This allows one to handle in a unified way approximation problems defined in terms of semifields other than max-plus algebra. Moreover, the new approach uses Puiseux polynomial and rational functions, which considerably extends the allowed classes of approximating functions and, hence, makes the approach more flexible and adaptable to a wider range of practical problems. The proposed solution of the rational approximation problem, which is based on a new efficient alternating algorithm for solving two-sided vector equations, and which seems to have no analogue in the literature, is one of the main advantages of the approach.

The rest of the paper proceeds as follows. In Section~\ref{S-DNPR}, we present the basic definitions and notation of tropical algebra, which are used throughout the paper. Section~\ref{S-BASVE} describes  the solution of one- and two-sided tropical linear vector equations. In Section~\ref{S-TDBAF}, we offer our main results, which provide solutions to discrete best approximation problems with tropical polynomial and rational approximating functions. Section~\ref{S-NE} includes numerical examples of solutions for approximation problems formulated in terms of max-plus algebra. In Section~\ref{S-C}, we give some concluding remarks and discuss further research directions.

\section{Definitions, Notation and Preliminary Results}
\label{S-DNPR}

We start with key definitions, notation and preliminary results that describe a formal framework for the development of solutions of approximation problems in the tropical algebra setting in what follows. Further information on tropical (idempotent) algebra and its applications can be found in a range of publications, including \cite{Kolokoltsov1997Idempotent,Golan2003Semirings,Heidergott2006Maxplus,Itenberg2007Tropical,Gondran2008Graphs,Butkovic2010Maxlinear,Maclagan2015Introduction}.

\subsection{Idempotent Semifield}
Let $\mathbb{X}$ be a non-empty set that is equipped with binary operations $\oplus$ (addition) and $\otimes$ (multiplication), and contains distinct elements $\mathbb{0}$ (zero) and $\mathbb{1}$ (unit). Assume that $(\mathbb{X},\oplus,\mathbb{0})$ is an idempotent commutative monoid, $(\mathbb{X}\setminus\{\mathbb{0}\},\otimes,\mathbb{1})$ is an Abelian group and multiplication $\otimes$ distributes over addition $\oplus$. We consider the algebraic system $(\mathbb{X},\oplus,\otimes,\mathbb{0},\mathbb{1})$, which is commonly referred to as a tropical (idempotent) semifield. 

The semifield has addition, which is idempotent: for each $x\in\mathbb{X}$ the equality $x\oplus x=x$ holds, and multiplication that is invertible: for each $x\ne\mathbb{0}$, there exists $x^{-1}$, such that $xx^{-1}=\mathbb{1}$ (here and hereafter the multiplication sign $\otimes$ is omitted, to save writing). Moreover, the semifield is considered complete, in the sense that the sums of infinite number of terms are well-defined and multiplication distributes over the infinite sums. 

The integer powers are defined in terms of multiplication $\otimes$, to represent repeated multiplication. It is assumed that the equation $x^{p}=a$ has a unique solution $x$ for any $a\in\mathbb{X}$ and integer $p>0$, which provides powers with rational exponents and, thus, turns the idempotent semifield into an algebraically closed (radicable) semifield.

Idempotent addition induces a partial order by the rule: $x\leq y$ if and only if $x\oplus y=y$. This partial order is assumed to extend to a total order, which makes the semifield linearly ordered. Since the semifield is complete, each non-empty set bounded from below has its greatest lower bound with respect to this order.

A typical example of the idempotent semifield under consideration is the real semifield $\mathbb{R}_{\max,+}=(\mathbb{R}\cup\{-\infty\},\max,+,-\infty,0)$, which is known as max-plus algebra. In this semifield, addition is defined as $\max$, multiplication as $+$, zero as $-\infty$ and unit as $0$. The power $x^{y}$ coincides with the arithmetic product $x\times y$ defined for any $x,y\in\mathbb{R}$. The inverse $x^{-1}$ of any $x\in\mathbb{R}$ corresponds to the usual opposite number $-x$. The order relation, which is associated with idempotent addition, agrees with the standard linear order on $\mathbb{R}$.

Another example is the semifield $\mathbb{R}_{\max}=(\mathbb{R}_{+},\max,\times,0,1)$, known as max-algebra, where $\mathbb{R}_{+}=\{x\in\mathbb{R}|x\geq0\}$. It has the operations $\oplus=\max$ and $\otimes=\times$ and the neutral elements $\mathbb{0}=0$ and $\mathbb{1}=1$. The notations of powers and inverses have the standard interpretation. The order induced by addition is the natural linear order on $\mathbb{R}_{+}$. 

In conclusion, we note that the algebraic expressions in what follows are thought of as defined in terms of an idempotent semifield unless otherwise indicated.

\subsection{Algebra of Matrices and Vectors}
Consider matrices and vectors over the idempotent semifield. The set of matrices of $M$ rows and $N$ columns with entries in $\mathbb{X}$ is denoted $\mathbb{X}^{M\times N}$. A matrix that has all entries equal to $\mathbb{0}$ is the zero matrix. If a matrix has no zero entries, it is called regular.

Addition, multiplication and scalar multiplication of matrices follow the standard entrywise rules, where the arithmetic addition and multiplication are replaced by $\oplus$ and $\otimes$. 

A matrix with a single column (row) is a column (row) vector. All vectors below are column vectors unless transposed. The set of column vectors that consist of $N$ elements is denoted by $\mathbb{X}^{N}$. A vector that has all elements equal to $\mathbb{0}$ is the zero vector denoted by $\bm{0}$. If a vector has no zero elements, it is called regular.

For any nonzero column vector $\bm{x}=(x_{i})$, the multiplicative conjugate is defined as the row vector $\bm{x}^{-}=(x_{i}^{-})$, with elements $x_{i}^{-}=x_{i}^{-1}$ if $x_{i}\ne\mathbb{0}$, and $x_{i}^{-}=\mathbb{0}$ otherwise.

\subsection{Tropical Vector Space}
Let $\bm{a}_{1},\ldots,\bm{a}_{N}\in\mathbb{X}^{M}$ be a set of vectors and $\bm{b}=x_{1}\bm{a}_{1}\oplus\cdots\oplus x_{N}\bm{a}_{N}$ be a linear combination of these vectors with coefficients $x_{1},\ldots,x_{N}\in\mathbb{X}$. The set of linear combinations $\mathcal{A}=\{x_{1}\bm{a}_{1}\oplus\cdots\oplus x_{N}\bm{a}_{N}|\ x_{1},\ldots,x_{N}\in\mathbb{X}\}$, which is closed under vector addition and scalar multiplication, is a tropical (finite-dimensional) vector space generated by $\bm{a}_{1},\ldots,\bm{a}_{N}$.

For any vector $\bm{a}=(a_{i})$, the support of $\bm{a}$ is given by $\mathop\mathrm{supp}(\bm{a})=\{i|\ a_{i}\ne\mathbb{0},\ 1\leq i\leq M\}$. Let us take two nonzero vectors $\bm{a}=(a_{i})$ and $\bm{b}=(b_{i})$, such that $\mathop\mathrm{supp}(\bm{a})=\mathop\mathrm{supp}(\bm{b})$, and let us define a distance function as follows:
\begin{equation}
d(\bm{a},\bm{b})
=
\bigoplus_{i\in\mathop\mathrm{supp}(\bm{a})}\left(a_{i}b_{i}^{-1}\oplus a_{i}^{-1} b_{i}\right)
=
\bm{b}^{-}\bm{a}
\oplus
\bm{a}^{-}\bm{b}.
\label{E-dab}
\end{equation}

When $\mathop\mathrm{supp}(\bm{a})\ne\mathop\mathrm{supp}(\bm{b})$, we set the function to a value that is greater than any $x\in\mathbb{X}$ (an undefined value) and write $d(\bm{a},\bm{b})=\infty$. If $\bm{a}=\bm{b}=\bm{0}$, then we assume $d(\bm{a},\bm{b})=\mathbb{1}$.
 
Suppose the function $d$ is defined in the framework of the semifield $\mathbb{R}_{\max,+}$ (max-plus algebra), where $\mathbb{1}=0$. In this context, it coincides for all $\bm{a},\bm{b}\in\mathbb{R}^{M}$ with the Chebyshev metric that is given in terms of the standard algebra as follows:
\begin{equation*}
d_{\infty}(\bm{a},\bm{b})
=
\max_{1\leq i\leq M}|a_{i}-b_{i}|
=
\max_{1\leq i\leq M}\max(a_{i}-b_{i},b_{i}-a_{i}).
\end{equation*}

In the case of the arbitrary idempotent semifield $\mathbb{X}$, we can interpret the distance function $d$ as a generalized metric that takes values in the subset $[\mathbb{1},\infty)$ of $\mathbb{X}$.

\subsection{Tropical Puiseux Polynomials and Rational Functions}

A tropical Puiseux monomial in one variable $x$ over $\mathbb{X}$ is a power $x^{p}$ with an exponent $p\in\mathrm{P}$, where $\mathrm{P}\subset\mathbb{Q}$ is a finite set of rationals. A tropical Puiseux polynomial is defined as a linear combination of monomials $x^{p}$ with given nonzero coefficients $\theta_{p}\in\mathbb{X}$ in the form
\begin{equation*}
P(x)
=
\bigoplus_{p\in\mathrm{P}}
\theta_{p}
x^{p},
\qquad
x\ne\mathbb{0}.
\end{equation*}

Since all coefficients satisfy the condition $\theta_{p}>\mathbb{0}$ and, hence, are tropically positive, the function $P$ can also be considered as a tropical posynomial.

Suppose that the polynomial $P$ consists of $N$ monomials, each given by its exponent $p$. Next, we assume that the monomials are ordered in increasing exponents $p$ and then consecutively numbered, starting from 1 to $N$. We use the numbers $i=1,\ldots,N$ to relabel the coefficients as $\theta_{i}$ and the exponents as $p_{i}$. Moreover, we can assume that $p_{i}\ne p_{j}$ for all $i\ne j$ (since, otherwise, one of the monomials $a_{i}x^{p_{i}}$ and $a_{j}x^{p_{i}}$ is always dominated by the other and can be omitted).

As a result, the tropical Puiseux polynomial can be written as
\begin{equation*}
P(x)
=
\bigoplus_{j=1}^{N}
\theta_{j}
x^{p_{j}},
\qquad
x\ne\mathbb{0},
\end{equation*}
where $\theta_{j}\ne\mathbb{0}$ and $p_{j}\in\mathbb{Q}$ for all $j=1,\ldots,N$, and $p_{1}<\cdots<p_{N}$.

We note that, in the context of the max-plus algebra (the semifield $\mathbb{R}_{\max,+}$), the polynomial is represented, in terms of the usual operations, as
\begin{equation*}
P(x)
=
\max_{1\leq j\leq N}
(p_{j}x+\theta_{j}),
\end{equation*} 
and, thus, defines a piecewise-linear convex function on $\mathbb{R}$. In the framework of $\mathbb{R}_{\max}$, the conventional form of the polynomial becomes
\begin{equation*}
P(x)
=
\max_{1\leq j\leq N}
\theta_{j}\times
x^{p_{j}},
\end{equation*} 
which in general specifies a nonlinear nonconvex function on $\mathbb{R}_{+}$.

Now, consider two polynomials:
\begin{equation*}
P(x)
=
\bigoplus_{j=1}^{N}
\theta_{j}
x^{p_{j}},
\qquad
Q(x)
=
\bigoplus_{k=1}^{L}
\sigma_{k}
x^{q_{k}},
\qquad
x\ne\mathbb{0},
\end{equation*}
and note that $Q(x)\ne\mathbb{0}$, since we assume that $\sigma_{k}\ne\mathbb{0}$ for all $k$.

A tropical rational function is the tropical quotient
\begin{equation*}
R(x)
=
\frac{P(x)}{Q(x)}.
\end{equation*}

If we interpret the rational function in terms of max-plus algebra, it is represented in the standard form as 
\begin{equation*}
R(x)
=
P(x)-Q(x)
=
\max_{1\leq j\leq N}
(p_{j}x+\theta_{j})
-
\max_{1\leq k\leq L}
(q_{k}x+\sigma_{k}),
\end{equation*}
which is a difference of convex functions (a DC-function). We observe that any arbitrary continuous function can be represented as the difference of two convex functions \cite{Tuy2009Dcprogramming}.

\section{Best Approximate Solution of Vector Equations}
\label{S-BASVE}

We are now concerned with the solution of linear vector equations, and we describe procedures to find best approximate solutions, in the sense of the distance function given by \eqref{E-dab}. We examine the equations 
\begin{equation*}
\bm{A}\bm{x}
=
\bm{b},
\qquad
\bm{A}\bm{x}
=
\bm{B}\bm{y},
\end{equation*}
where $\bm{A}$ and $\bm{B}$ are given matrices, $\bm{b}$ is a given vector and $\bm{x}$ and $\bm{y}$ are unknown vectors.

The first equation, which has an unknown vector on one side, is normally referred to as a one-sided equation. The second equation, with unknowns on both sides, is two-sided.

To solve the equations, we first evaluate distances between vectors in tropical vector spaces, and we then apply the result, to derive the best approximate solutions.

\subsection{Solution of One-Sided Equation}

Consider a tropical vector space $\mathcal{A}$ generated by nonzero vectors $\bm{a}_{1},\ldots,\bm{a}_{N}\in\mathbb{X}^{M}$. Any vector $\bm{a}\in\mathcal{A}$ can be written as a linear combination $\bm{a}=x_{1}\bm{a}_{1}\oplus\cdots\oplus x_{N}\bm{a}_{N}$, where $x_{1},\ldots,x_{N}\in\mathbb{X}$, and can then be represented in the form of the matrix-vector product as $\bm{a}=\bm{A}\bm{x}$, with the matrix $\bm{A}=(\bm{a}_{1},\ldots,\bm{a}_{N})$ and vector $\bm{x}=(x_{1},\ldots,x_{N})^{T}$.

We examine the distance from a vector $\bm{b}$ to the vector space $\mathcal{A}$, which is given by
\begin{equation*}
d(\mathcal{A},\bm{b})
=
\min_{\bm{a}\in\mathcal{A}}d(\bm{a},\bm{b})
=
\min_{\bm{x}\in\mathbb{X}^{n}}d(\bm{A}\bm{x},\bm{b}).
\end{equation*}

In the case when vector $\bm{b}$ is regular, we can replace the minimum over all $\bm{x}\in\mathbb{X}$ on the right-hand side by the minimum taken over regular $\bm{x}$ (see, e.g., \cite{Krivulin2009Solution}) and write
\begin{equation*}
d(\mathcal{A},\bm{b})
=
\min_{\bm{x}>\bm{0}}d(\bm{A}\bm{x},\bm{b})
=
\min_{\bm{x}>\bm{0}}(\bm{b}^{-}\bm{A}\bm{x}\oplus(\bm{A}\bm{x})^{-}\bm{b}).
\end{equation*}

It follows from the definition of $d$ that the equality $d(\mathcal{A},\bm{b})=\mathbb{1}$ indicates the condition $\bm{b}\in\mathcal{A}$, whereas the inequality $d(\mathcal{A},\bm{b})>\mathbb{1}$ corresponds to the condition $\bm{b}\not\in\mathcal{A}$.

The next result is obtained in \cite{Krivulin2009Solution} (see also \cite{Krivulin2013Solution-linear}):
\begin{lemma}
\label{L-mindAxb}
Let $\bm{A}$ be a regular matrix, $\bm{b}$ a regular vector and denote $\Delta=(\bm{A}(\bm{b}^{-}\bm{A})^{-})^{-}\bm{b}$. Then,
\begin{equation*}
\min_{\bm{x}>\bm{0}}\ d(\bm{A}\bm{x},\bm{b})
=
\sqrt{\Delta},
\end{equation*}
where the minimum is achieved at $\bm{x}=\sqrt{\Delta}(\bm{b}^{-}\bm{A})^{-}$.
\end{lemma}

The above lemma demonstrates that the distance from a vector $\bm{b}$ to the vector space $\mathcal{A}$ generated by the columns of a matrix $\bm{A}$ is calculated as $d(\mathcal{A},\bm{b})=\sqrt{\Delta}$. Moreover, the vector given by $\bm{y}=\sqrt{\Delta}\bm{A}(\bm{b}^{-}\bm{A})^{-}$ is a vector of $\mathcal{A}$, which is closest to the vector $\bm{b}$.  

Given a matrix $\bm{A}\in\mathbb{X}^{M\times N}$ and vector $\bm{b}\in\mathbb{X}^{M}$, consider the problem of finding regular vectors $\bm{x}\in\mathbb{X}^{N}$ that solve the one-sided equation
\begin{equation}
\bm{A}\bm{x}
=
\bm{b}.
\label{E-Axeqb}
\end{equation}

Observing that the equation can have no solution (be inconsistent), we turn to finding a best approximate solution, in the sense of the metric $d$. Vector $\bm{x}_{\ast}$ is a best approximate solution of Equation \eqref{E-Axeqb} if it satisfies for all vectors $\bm{x}$ the inequality
\begin{equation*}
d(\bm{A}\bm{x}_{\ast},\bm{b})
\leq
d(\bm{A}\bm{x},\bm{b}).
\end{equation*}

We note that vector $\bm{A}\bm{x}_{\ast}$ can be considered as a projection of $\bm{b}$ onto the space $\mathcal{A}$.

A solution to Equation \eqref{E-Axeqb} can be obtained, based on results in \cite{Krivulin2009Solution,Krivulin2013Solution-linear}, as follows.
\begin{theorem}
\label{T-Axeqb}
Let $\bm{A}$ be a regular matrix, $\bm{b}$ a regular vector and denote $\Delta=(\bm{A}(\bm{b}^{-}\bm{A})^{-})^{-}\bm{b}$. Then:
\begin{enumerate}
\item
The best approximate solution of Equation \eqref{E-Axeqb}, in the sense of the distance function $d$, is \mbox{given by}
\begin{equation*}
\bm{x}_{\ast}
=
\sqrt{\Delta}(\bm{b}^{-}\bm{A})^{-}.
\end{equation*}
\item
If $\Delta=\mathbb{1}$, the equation has exact solutions; the vector $\bm{x}_{\ast}=(\bm{b}^{-}\bm{A})^{-}$ is the maximal solution.
\end{enumerate}
\end{theorem}

Since the value of $\sqrt{\Delta}$ shows the minimum deviation between both sides of \eqref{E-Axeqb}, measured on the scale of the function $d$, it also specifies the approximation error.

\subsection{Solution of Two-Sided Equation}

Let $\mathcal{A}$ be a vector space generated by vectors $\bm{a}_{1},\ldots,\bm{a}_{N}\in\mathbb{X}^{M}$, and let $\mathcal{B}$ be a vector space generated by $\bm{b}_{1},\ldots,\bm{b}_{L}\in\mathbb{X}^{M}$. With the matrices $\bm{A}=(\bm{a}_{1},\ldots,\bm{a}_{N})$ and $\bm{B}=(\bm{b}_{1},\ldots,\bm{b}_{L})$, we define the distance between the spaces as
\begin{equation*}
d(\mathcal{A},\mathcal{B})
=
\min_{\bm{a}\in\mathcal{A},\bm{b}\in\mathcal{B}}d(\bm{a},\bm{b})
=
\min_{\bm{x}>\bm{0},\bm{y}>\bm{0}}d(\bm{A}\bm{x},\bm{B}\bm{y})
=
\min_{\bm{x}>\bm{0}}((\bm{B}\bm{y})^{-}\bm{A}\bm{x}\oplus(\bm{A}\bm{x})^{-}\bm{B}\bm{y}).
\end{equation*}

The equality $d(\mathcal{A},\mathcal{B})=\mathbb{1}$ indicates that the vector spaces $\mathcal{A}$ and $\mathcal{B}$ have non-empty intersection. If the inequality $d(\mathcal{A},\mathcal{B})>\mathbb{1}$ holds, then these spaces have no common point. In this case, the value of the function $d$ shows the minimum distance between the spaces.

Suppose that $\bm{A}\in\mathbb{X}^{M\times N}$ and $\bm{B}\in\mathbb{X}^{M\times L}$ are given regular matrices. Consider the problem of finding regular vectors $\bm{x}\in\mathbb{X}^{N}$ and $\bm{y}\in\mathbb{X}^{L}$ that solve the two-sided equation,
\begin{equation}
\bm{A}\bm{x}
=
\bm{B}\bm{y}.
\label{E-AxeqBy}
\end{equation}

A pair of vectors $(\bm{x}_{\ast},\bm{y}_{\ast})$ is a best approximate solution of Equation \eqref{E-AxeqBy} if, for all vectors $\bm{x}$ and $\bm{y}$, the following inequality holds:
\begin{equation*}
d(\bm{A}\bm{x}_{\ast},\bm{B}\bm{y}_{\ast})
\leq
d(\bm{A}\bm{x},\bm{B}\bm{y}).
\end{equation*}

We note that the problem of finding a best approximate solution of Equation \eqref{E-AxeqBy} is close to the best proximity point problems introduced as an extension of the classical fixed point problem. For further details on the related theory and references, one can see \cite{Basha2011Best}.

To obtain a best approximate solution, we apply the procedure proposed in \cite{Krivulin2023Solution}, which constructs a sequence of vectors from the spaces $\mathcal{A}$ and $\mathcal{B}$ generated by the columns of the matrices $\bm{A}$ and $\bm{B}$. The vectors are taken alternatively from both spaces, in such a way that after selecting a vector in one space the next vector is found in the other space, to minimize the distance to the former vector. The vectors in the spaces $\mathcal{A}$ and $\mathcal{B}$ are given by coefficients in their decompositions as linear combinations of columns in $\bm{A}$ and $\bm{B}$.

Let $\bm{x}_{0}\in\mathbb{X}^{N}$ be a regular vector and $\bm{a}_{0}=\bm{A}\bm{x}_{0}\in\mathcal{A}$. We apply Theorem~\ref{T-Axeqb}, to find the minimum distance from the vector $\bm{a}_{0}$ to the vectors in $\mathcal{B}$ as follows:
\begin{equation*}
d(\bm{a}_{0},\mathcal{B})
=
\sqrt{\Delta_{0}},
\qquad
\Delta_{0}
=
(\bm{B}((\bm{A}\bm{x}_{0})^{-}\bm{B})^{-})^{-}\bm{A}\bm{x}_{0}.
\end{equation*}

This minimum distance is provided by a vector $\bm{b}_{1}\in\mathcal{B}$ that is calculated as
\begin{equation*}
\bm{b}_{1}
=
\bm{B}\bm{y}_{1},
\qquad
\bm{y}_{1}
=
\sqrt{\Delta_{0}}((\bm{A}\bm{x}_{0})^{-}\bm{B})^{-}.
\end{equation*}

We evaluate the minimum distance from $\bm{b}_{1}$ to the vectors in $\mathcal{A}$, which is given by 
\begin{equation*}
d(\bm{b}_{1},\mathcal{A})
=
\sqrt{\Delta_{1}},
\qquad
\Delta_{1}
=
(\bm{A}((\bm{B}\bm{y}_{1})^{-}\bm{A})^{-})^{-}\bm{B}\bm{y}_{1},
\end{equation*}
and attained at a vector $\bm{a}_{2}\in\mathcal{A}$, such that
\begin{equation*}
\bm{a}_{2}
=
\bm{A}\bm{x}_{2},
\qquad
\bm{x}_{2}
=
\sqrt{\Delta_{1}}((\bm{B}\bm{y}_{1})^{-}\bm{A})^{-}.
\end{equation*}

Furthermore, we evaluate the distance $d(\bm{a}_{2},\mathcal{B})$, by calculating $\Delta_{2}$, and we then use this distance to obtain vectors $\bm{y}_{3}$ and $\bm{b}_{3}$. Next, we evaluate $\Delta_{3}$ and find vectors $\bm{x}_{4}$ and $\bm{a}_{4}$.

We continue the alternating procedure, to form a sequence of vectors $\bm{a}_{0},\bm{b}_{1},\bm{a}_{2},\bm{b}_{3},\bm{a}_{4},\ldots$ that are taken alternatively from the spaces $\mathcal{A}$ and $\mathcal{B}$, to minimize the distance between successive vectors. A sequence of pairs $(\bm{x}_{0},\bm{y}_{1}),(\bm{x}_{2},\bm{y}_{3}),\ldots$ is also generated, which provides successive approximations to the best approximate solution of Equation \eqref{E-AxeqBy}.

Consider the sequence $\Delta_{0},\Delta_{1},\Delta_{2},\ldots$, which is bounded from below, since $\Delta_{i}\geq\mathbb{1}$ for all $i=0,1,2,\ldots$. One can verify (see \cite{Krivulin2023Solution}) that $\Delta_{i+1}\leq\Delta_{i}$ and, hence, the sequence is non-increasing. As a result, this sequence converges to a limit $\Delta_{\ast}\geq\mathbb{1}$, defined as the greatest lower bound of its elements. Observing that the underlying semifield is complete, we conclude that this limit exists and coincides with the minimum element of the sequence.

We note that each element of the sequence specifies the squared distance between a vector of one of the spaces $\mathcal{A}$ and $\mathcal{B}$ and the nearest vector in the other space. If the equality $\Delta_{i}=\mathbb{1}$ is satisfied for some $i$, then the spaces have non-empty intersection, and \mbox{Equation \eqref{E-AxeqBy}} has exact regular solutions. In the case of even $i$, the intersection contains the vector $\bm{a}_{i}=\bm{A}\bm{x}_{i}$, whereas the pair of vectors $(\bm{x}_{i},\bm{y}_{i+1})$ is a solution. If $i$ is odd, the intersection contains $\bm{b}_{i}=\bm{B}\bm{y}_{i}$, whereas the pair $(\bm{x}_{i+1},\bm{y}_{i})$ is a solution.

If the inequality $\Delta_{\ast}>\mathbb{1}$ holds, then the spaces $\mathcal{A}$ and $\mathcal{B}$ do not intersect, and \mbox{Equation \eqref{E-AxeqBy}} does not have exact solutions. In such an event, the procedure stops as soon as a repeated element occurs in any of the sequences $\bm{x}_{0},\bm{x}_{2},\ldots$ or $\bm{y}_{1},\bm{y}_{3},\ldots$.

The above alternating procedure can be represented as the following algorithm.

\begin{algorithm}

\textit{Solution of two-sided equation $\bm{A}\bm{x}=\bm{B}\bm{x}$}
\begin{enumerate}
\item
Input regular matrices $\bm{A},\bm{B}$ and regular vector $\bm{x}_{0}$; set $i=0$.
\item\label{Loop}
Calculate
\begin{equation*}
\Delta_{i}
=
(\bm{B}((\bm{A}\bm{x}_{i})^{-}\bm{B})^{-})^{-}\bm{A}\bm{x}_{i},
\qquad
\bm{y}_{i+1}
=
\sqrt{\Delta_{i}}((\bm{A}\bm{x}_{i})^{-}\bm{B})^{-}.
\end{equation*}
\item
If $\Delta_{i}=\mathbb{1}$ or $\bm{y}_{i+1}=\bm{y}_{j}$ for some $j<i$, then set
\begin{equation*}
\Delta_{\ast}
=
\Delta_{i},
\qquad
\bm{x}_{\ast}
=
\bm{x}_{i},
\qquad
\bm{y}_{\ast}
=
\bm{y}_{i+1},
\end{equation*}
and stop; otherwise, set $i=i+1$.
\item
Calculate
\begin{equation*}
\Delta_{i}
=
(\bm{A}((\bm{B}\bm{y}_{i})^{-}\bm{A})^{-})^{-}\bm{B}\bm{y}_{i},
\qquad
\bm{x}_{i+1}
=
\sqrt{\Delta_{i}}((\bm{B}\bm{y}_{i})^{-}\bm{A})^{-}.
\end{equation*}
\item
If $\Delta_{i}=\mathbb{1}$ or $\bm{x}_{i+1}=\bm{x}_{j}$ for some $j<i$, then set
\begin{equation*}
\Delta_{\ast}
=
\Delta_{i},
\qquad
\bm{x}_{\ast}
=
\bm{x}_{i+1},
\qquad
\bm{y}_{\ast}
=
\bm{y}_{i},
\end{equation*}
and stop; otherwise, set $i=i+1$.
\item
Go to step \ref{Loop}.
\end{enumerate}
 \label{A-AxeqBy}
\end{algorithm}

Suppose that upon completion of the algorithm, we obtain $\Delta_{\ast}=\mathbb{1}$. This result means that Equation \eqref{E-AxeqBy} has exact regular solutions, which include the obtained pair of vectors $(\bm{x}_{\ast},\bm{y}_{\ast})$. If we have $\Delta_{\ast}>\mathbb{1}$, then the equation has no exact solution. The value of $\Delta_{\ast}$ specifies the approximation error attained at a best approximate solution given by $(\bm{x}_{\ast},\bm{y}_{\ast})$.

\section{Tropical Discrete Best Approximation of Functions}
\label{S-TDBAF}

We are now in a position to propose new algebraic techniques to solve the data approximation problems of finding an unknown but definite function $y=f(x)$ from finitely many samples $(x_{i},y_{i})$ in the tropical algebra setting. The purpose is to determine a functional model that can effectively approximate the underlying relation between $x_{i}$ and $y_{i}$.

As approximating functions, the techniques use tropical polynomial and rational functions. The approximation problems are reduced to evaluating the best approximate solutions of vector equations obtained from the sample data.

\subsection{Tropical Polynomial Approximation}

Suppose that samples $x_{i},y_{i}\in\mathbb{X}$ are given for $i=1,\ldots,M$, where $y_{i}$ is assumed to be an output obtained from an input $x_{i}$ by some function $f:\mathbb{X}\rightarrow\mathbb{X}$. Let us consider the problem of approximating the function $f$ by a polynomial, using these sample data. 

For each $j=1,\ldots,N$, let $p_{j}\in\mathbb{Q}$ denote a fixed number, and $\theta_{j}\in\mathbb{X}$ denote an unknown parameter. Define the approximating polynomial as
\begin{equation*}
P(x)
=
\bigoplus_{j=1}^{N}\theta_{j}x^{p_{j}}
=
\theta_{1}x^{p_{1}}
\oplus\cdots\oplus
\theta_{N}x^{p_{N}}.
\end{equation*}

The problem consists in the determination of the parameters that make the equations
\begin{equation*}
\theta_{1}x_{i}^{p_{1}}
\oplus\cdots\oplus
\theta_{N}x_{i}^{p_{N}}
=
y_{i}
\qquad
i=1,\ldots,M,
\end{equation*}
hold exactly or approximately, by minimizing the deviation between both sides. 

We introduce the matrix-vector notation,
\begin{equation*}
\bm{X}
=
\left(
\begin{array}{ccc}
x_{1}^{p_{1}} & \ldots & x_{1}^{p_{N}}
\\
\vdots & & \vdots
\\
x_{M}^{p_{1}} & \ldots & x_{M}^{p_{N}}
\end{array}
\right),
\qquad
\bm{y}
=
\left(
\begin{array}{c}
y_{1},
\\
\vdots
\\
y_{M}
\end{array}
\right),
\qquad
\bm{\theta}
=
\left(
\begin{array}{c}
\theta_{1}
\\
\vdots
\\
\theta_{N}
\end{array}
\right),
\end{equation*}
and we combine the above scalar equations into the one-sided vector equation,
\begin{equation*}
\bm{X}\bm{\theta}
=
\bm{y},
\end{equation*}
where $\bm{X}$ and $\bm{y}$ are a known matrix and vector, and $\bm{\theta}$ is an unknown vector of parameters.

To make the equation fit the sample data, we need to find a best approximate solution, in the sense of a minimum distance between the vectors $\bm{X}\bm{\theta}$ and $\bm{y}$. With the distance function $d$, the problem is to obtain a vector $\bm{\theta}_{\ast}=(\theta_{1}^{\ast},\ldots,\theta_{N}^{\ast})^{T}$ that is given by
\begin{equation*}
\bm{\theta}_{\ast}
=
\arg\min_{\bm{\theta}}
d(\bm{X}\bm{\theta},\bm{y}).
\end{equation*}
  
We solve the problem by applying Theorem~\ref{T-Axeqb}, which provides a best approximation solution in terms of the metric $d$. We replace the matrix $\bm{A}$ by $\bm{X}$, and the vectors $\bm{b}$ and $\bm{x}$ by $\bm{y}$ and $\bm{\theta}$. First, we calculate the squared best approximation error:  
\begin{equation*}
\Delta_{\ast}
=
(\bm{X}(\bm{y}^{-}\bm{X})^{-})^{-}\bm{y}.
\end{equation*}

Next, we obtain the vector of parameters that achieves this error:
\begin{equation*}
\bm{\theta}_{\ast}
=
\sqrt{\Delta_{\ast}}(\bm{y}^{-}\bm{X})^{-}.
\end{equation*}

The best approximating polynomial is then given by
\begin{equation*}
P_{\ast}(x)
=
\theta_{1}^{\ast}x^{p_{1}}
\oplus\cdots\oplus
\theta_{N}^{\ast}x^{p_{N}}.
\end{equation*}

\subsection{Tropical Rational Approximation}

As before, let $x_{i},y_{i}\in\mathbb{X}$ be given samples for $i=1,\ldots,M$. Furthermore, let $p_{j},q_{k}\in\mathbb{Q}$ be fixed numbers and $\theta_{j},\sigma_{k}\in\mathbb{X}$ be unknown parameters for $j=1,\ldots,N$ and $k=1,\ldots,L$.

Denote numerator and denominator polynomials by
\begin{equation*}
P(x)
=
\bigoplus_{j=1}^{N}\theta_{j}x^{p_{j}},
\qquad
Q(x)
=
\bigoplus_{k=1}^{L}\sigma_{k}x^{q_{k}},
\end{equation*}
and define an approximating rational function as
\begin{equation*}
R(x)
=
\frac{P(x)}{Q(x)}
=
\frac{\theta_{1}x^{p_{1}}\oplus\cdots\oplus\theta_{N}x^{p_{N}}}{\sigma_{1}x^{q_{1}}\oplus\cdots\oplus\sigma_{L}x^{q_{L}}}.
\end{equation*}

We now consider the problem of evaluating the parameters that achieve the best approximation of the equations
\begin{equation}
\frac{\theta_{1}x_{i}^{p_{1}}\oplus\cdots\oplus\theta_{N}x_{i}^{p_{N}}}{\sigma_{1}x_{i}^{q_{1}}\oplus\cdots\oplus\sigma_{L}x_{i}^{q_{L}}}
=
y_{i},
\qquad
i=1,\ldots,M.
\label{E-XthetaZsigmaeqy}
\end{equation}
 
We rewrite the equations as
\begin{equation*}
\theta_{1}x_{i}^{p_{1}}\oplus\cdots\oplus\theta_{N}x_{i}^{p_{N}}
=
y_{i}(\sigma_{1}x_{i}^{q_{1}}\oplus\cdots\oplus\sigma_{L}x_{i}^{q_{L}}),
\qquad
i=1,\ldots,M.
\end{equation*}

To represent the problem in vector form, we use the matrix-vector notation
\begin{gather*}
\bm{X}
=
\left(
\begin{array}{ccc}
x_{1}^{p_{1}} & \ldots & x_{1}^{p_{N}}
\\
\vdots & & \vdots
\\
x_{M}^{p_{1}} & \ldots & x_{M}^{p_{N}}
\end{array}
\right),
\qquad
\bm{Y}
=
\left(
\begin{array}{ccc}
y_{1} &  & \mathbb{0}
\\
& \ddots &
\\
\mathbb{0} & & y_{M}
\end{array}
\right),
\\
\bm{Z}
=
\left(
\begin{array}{ccc}
x_{1}^{q_{1}} & \ldots & x_{1}^{q_{L}}
\\
\vdots & & \vdots
\\
x_{M}^{q_{1}} & \ldots & x_{M}^{q_{L}}
\end{array}
\right),
\qquad
\bm{\theta}
=
\left(
\begin{array}{c}
\theta_{1}
\\
\vdots
\\
\theta_{N}
\end{array}
\right),
\qquad
\bm{\sigma}
=
\left(
\begin{array}{c}
\sigma_{1}
\\
\vdots
\\
\sigma_{L}
\end{array}
\right),
\end{gather*}
where $\bm{X}$, $\bm{Y}$ and $\bm{Z}$ are known data matrices, and where $\bm{\theta}$ and $\bm{\sigma}$ are unknown vectors of parameters. 

With this notation, the scalar equations take the form of the two-sided vector equation
\begin{equation*}
\bm{X}\bm{\theta}
=
\bm{Y}\bm{Z}\bm{\sigma}.
\end{equation*}

Let us verify that the approximation error for the last vector equation coincides with that of the system \eqref{E-XthetaZsigmaeqy}. We define the vector $\bm{w}=(w_{i})$ with elements $w_{i}=(\bm{X}\bm{\theta})_{i}(\bm{Z}\bm{\sigma})_{i}^{-1}$, and we examine the distance $d(\bm{w},\bm{y})$ that represents the approximation error for \eqref{E-XthetaZsigmaeqy}. By using the scalar representation at \eqref{E-dab}, we immediately obtain
\begin{equation*}
d(\bm{w},\bm{y})
=
\bigoplus_{i=1}^{M}((\bm{X}\bm{\theta})_{i}(\bm{Z}\bm{\sigma})_{i}^{-1}y_{i}^{-1}\oplus(\bm{X}\bm{\theta})_{i}^{-1}(\bm{Z}\bm{\sigma})_{i}y_{i})
=
d(\bm{X}\bm{\theta},\bm{Y}\bm{Z}\bm{\sigma}).
\end{equation*}

We now seek for parameters that provide a minimum distance between $\bm{X}\bm{\theta}$ and $\bm{Y}\bm{Z}\bm{\sigma}$. We need to find vectors $\bm{\theta}_{\ast}=(\theta_{1}^{\ast},\ldots,\theta_{N}^{\ast})^{T}$ and $\bm{\sigma}_{\ast}=(\sigma_{1}^{\ast},\ldots,\sigma_{L}^{\ast})^{T}$, such that
\begin{equation*}
(\bm{\theta}_{\ast},\bm{\sigma}_{\ast})
=
\arg\min_{\bm{\theta},\bm{\sigma}}
d(\bm{X}\bm{\theta},\bm{Y}\bm{Z}\bm{\sigma}).
\end{equation*}
 
To solve the problem, we apply Algorithm~\ref{A-AxeqBy}, where the matrices $\bm{A}$ and $\bm{B}$ are replaced by $\bm{X}$ and $\bm{Y}\bm{Z}$, and the vectors $\bm{x}$ and $\bm{y}$ by $\bm{\theta}$ and $\bm{\sigma}$. The algorithm yields the squared minimum distance $\Delta_{\ast}$ and vectors $\bm{\theta}^{\ast}$ and $\bm{\sigma}_{\ast}$ that attain the minimum.\vspace{+10pt}

The best approximate rational function is given by
\begin{equation*}
R_{\ast}(x)
=
\frac{\theta_{1}^{\ast}x^{p_{1}}
\oplus\cdots\oplus
\theta_{N}^{\ast}x^{p_{N}}}
{\sigma_{1}^{\ast}x^{q_{1}}
\oplus\cdots\oplus
\sigma_{L}^{\ast}x^{q_{L}}}.
\end{equation*}

To conclude this section, we note that the proposed approximation techniques assume that the number of terms (monomials) and their related degrees in the approximating polynomials $P(x)$ and $Q(x)$ are fixed in advance. For some approximation problems, these characteristics of polynomials can be determined based on additional information about the nature and properties of the phenomena under study. However, in many cases there is no prior information on the characteristics, which makes the selection of approximating polynomials a nontrivial challenge that needs to be addressed in future research.

As a first simple but reasonable step to meet the challenge, a Monte Carlo random sampling technique can be implemented, to search for the class of optimal polynomials. We use a random sampling technique to determine the degrees of monomials when solving illustrative example problems in the next section.

\section{Numerical Examples}
\label{S-NE}

In this section, we demonstrate how tropical approximation techniques can be used to solve discrete approximation problems with sample data evaluated from simple convex and nonconvex functions. To provide transparent graphical illustrations, we offer examples, in terms of the semifield $\mathbb{R}_{\max,+}$ (max-plus algebra), where both polynomials and rational functions can be represented in Cartesian coordinates as piecewise linear functions.

We consider approximation problems in which approximating polynomials have a fixed number of monomials, while the degrees of these monomials are not given in advance. We propose a two-level solution procedure that combines direct random search, to select a class of approximating polynomials or rational functions at the upper level, with best approximation using approximants from the selected class at the lower level.

To solve polynomial approximation problems, given a number of monomials in the polynomial, we implement uniform random sampling over a feasible set of degrees. Each sample sets the degree of each monomial to a fixed value and, hence, we define a class of polynomials with specified degrees and arbitrary coefficients. Then, we apply Theorem~\ref{T-Axeqb}, to find coefficients of a best approximate polynomial in the class, and we evaluate the related approximation error, to allow the search of a sample that produces the minimum error.

The solution of rational approximation problems includes independent random sampling for both polynomials in the numerator and denominator of the approximating rational function. The coefficients of the polynomials and corresponding approximation error are obtained by solving the approximation problem, with the use of Algorithm~\ref{A-AxeqBy}.

To reduce the number of classes of polynomials, to examine and simplify the random search in numerical experiments, we consider only polynomials with integer degrees.

\subsection{Discrete Approximation of Convex Function}

We consider a function defined in terms of standard algebra as
\begin{equation*}
f(x)
=
x^{2}-3x^{1/3}+5/2,
\qquad
0\leq x\leq 2.
\end{equation*}

We suppose that there are $M=21$ samples $(x_{i},y_{i})$ available, obtained with $f(x)$ on the interval $[0,2]$ at $x_{i}=(i-1)/10$, where $i=1,\ldots,M$. The problem is to find an approximate tropical polynomial that best describes the relationship between $x_{i}$ and $y_{i}$ for all $i$.

Observing that the function $f(x)$ is convex on the interval considered, we start with the approximation of $f(x)$ by polynomials. We first consider polynomials of $N=5$ terms,
\begin{equation*}
P(x)
=
\theta_{1}x^{p_{1}}\oplus\cdots\oplus\theta_{5}x^{p_{5}},
\end{equation*}
where the vector of degrees $\bm{p}=(p_{1},\ldots,p_{5})^{T}$ is not fixed in advance and, thus, its elements are produced by random sampling from the discrete uniform distribution over $[-15,5]$.

For each obtained sample vector of degrees, we evaluate the vector of coefficients $\bm{\theta}=(\theta_{1},\ldots,\theta_{5})^{T}$, to achieve the minimum approximation error corresponding to the sample. The combination of sample degrees and their related coefficients that attains the minimum error over all samples is taken, to form a best approximate polynomial. 

After examining 10,000 sample vectors of exponents, we arrive at the minimum squared error $\Delta_{\ast}=0.1360$. We obtain the following vectors of coefficients and degrees:
\begin{gather*}
\bm{\theta}_{\ast}^{T}
=
\left(
\begin{array}{ccccc}
2.5680 & 0.9176 & -0.4320 & -1.6281 & -3.2413
\end{array}
\right),
\\
\bm{p}_{\ast}^{T}
=
\left(
\begin{array}{ccccc}
-14 & -1 & 1 & 2 & 3
\end{array}
\right).
\end{gather*}

The obtained polynomial represented in the conventional form is given by
\begin{multline*}
P_{\ast}(x)
=
\max(2.5680-15x,
0.9176-x, 
-0.4320+x,    
-1.6281+2x,    
-3.2413+3x).
\end{multline*}

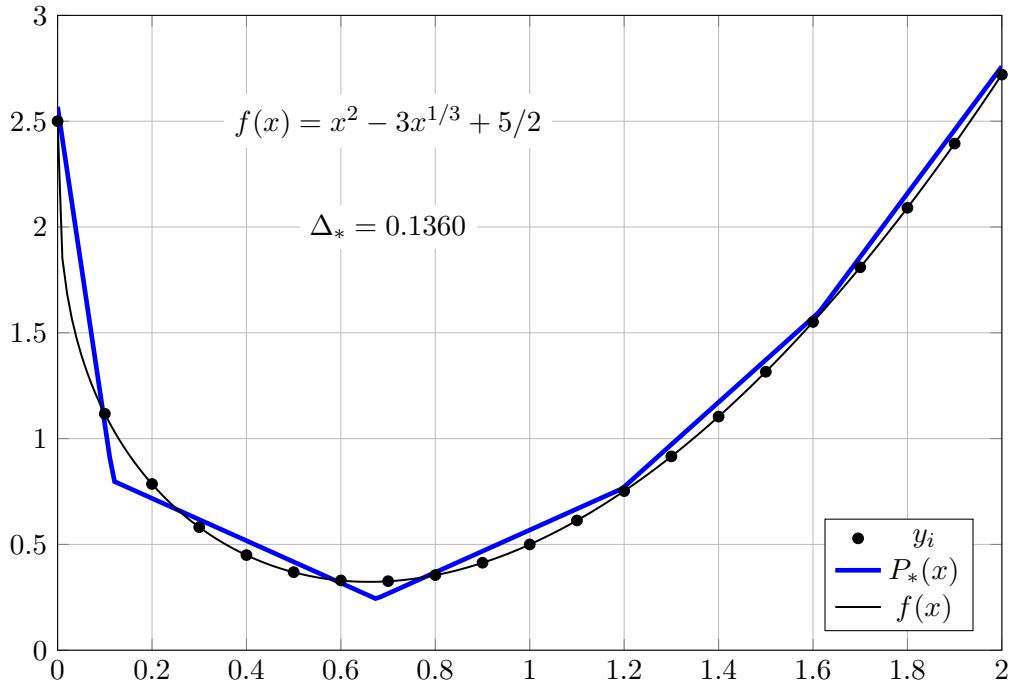
\begin{figure}[h]
\begin{tikzpicture}

\begin{axis}[legend pos=south east,
width=14cm,
height=10cm,
grid=major,
ymin=0,
ymax=3,
xmin=0,
xmax=2,
]

\addplot[
only marks,
]
coordinates {
(0, 2.5000)
(0.1000, 1.1175)
(0.2000, 0.7856)
(0.3000, 0.5817)
(0.4000, 0.4496)
(0.5000, 0.3689)
(0.6000, 0.3297)
(0.7000, 0.3263)
(0.8000, 0.3550)
(0.9000, 0.4135)
(1.0000, 0.5000)
(1.1000, 0.6132)
(1.2000, 0.7520)
(1.3000, 0.9158)
(1.4000, 1.1039)
(1.5000, 1.3159)
(1.6000, 1.5512)
(1.7000, 1.8096)
(1.8000, 2.0907)
(1.9000, 2.3943)
(2.0000, 2.7202)
};

\addlegendentry{$y_{i}$}

\addplot[
blue,
samples=200,
line width=1.75pt,
domain=0.0:2.0,
y domain=0:10,
]{max(
2.568004+x*(-15),
0.917585+x*(-1),
(-0.431996)+x*1,
(-1.628063)+x*2,
(-3.241317)+x*3
)};

\addlegendentry{$P_{\ast}(x)$}

\addplot[
samples=200,
black,
line width=0.75pt,
domain=0.0:2.0,
y domain=0:10,
]{x^2-3*x^(1/3)+5/2};

\addlegendentry{$f(x)$}

\node[style={fill=white}] at (axis cs: 0.7,2.5) {$f(x)=x^{2}-3x^{1/3}+5/2$};
\node[style={fill=white}] at (axis cs: 0.7,2.0) {$\Delta_{\ast}=0.1360$};

\end{axis}

\end{tikzpicture}
\caption{Approximation of the function $f(x)$ by a polynomial $P_{\ast}(x)$ with $N=5$ terms.}
\label{F-AFPN5}
\end{figure}

A graphical illustration of the solution is presented in Figure~\ref{F-AFPN5}.

Next, we use polynomials with $N=7$ terms, in the form
\begin{equation*}
P(x)
=
\theta_{1}x^{p_{1}}\oplus\cdots\oplus\theta_{7}x^{p_{7}}.
\end{equation*}

Approximation by these polynomials, combined with a random search based on 10,000 samples, yields the minimum $\Delta_{\ast}=0.0481$ and the vectors 
\begin{gather*}
\bm{\theta}_{\ast}^{T}
=
\left(
\begin{array}{ccccccc}
2.5240 & 1.4096 & 0.8736 & 0.3503 & -0.4760 & -1.6720 & -3.2853
\end{array}
\right),
\\
\bm{p}_{\ast}^{T}
=
\left(
\begin{array}{ccccccc}
-15 & -3 & -1 & 0 & 1 & 2 & 3
\end{array}
\right).
\end{gather*}

\begin{figure}[h]
\begin{tikzpicture}

\begin{axis}[legend pos=south east,
width=14cm,
height=10cm,
grid=major,
ymin=0,
ymax=3,
xmin=0,
xmax=2,
]

\addplot[
only marks,
]
coordinates {
(0, 2.5000)
(0.1000, 1.1175)
(0.2000, 0.7856)
(0.3000, 0.5817)
(0.4000, 0.4496)
(0.5000, 0.3689)
(0.6000, 0.3297)
(0.7000, 0.3263)
(0.8000, 0.3550)
(0.9000, 0.4135)
(1.0000, 0.5000)
(1.1000, 0.6132)
(1.2000, 0.7520)
(1.3000, 0.9158)
(1.4000, 1.1039)
(1.5000, 1.3159)
(1.6000, 1.5512)
(1.7000, 1.8096)
(1.8000, 2.0907)
(1.9000, 2.3943)
(2.0000, 2.7202)
};

\addlegendentry{$y_{i}$}

\addplot[
blue,
samples=200,
line width=1.75pt,
domain=0.0:2.0,
y domain=0:10,
]{max(
2.524046+x*(-15),
1.409635+x*(-3),
0.873627+x*(-1),
0.350334+x*0,
(-0.475954)+x*1,
(-1.672021)+x*2,
(-3.285276)+x*3
)};

\addlegendentry{$P_{\ast}(x)$}

\addplot[
samples=200,
black,
line width=0.75pt,
domain=0.0:2.0,
y domain=0:10,
]{x^2-3*x^(1/3)+5/2};

\addlegendentry{$f(x)$}

\node[style={fill=white}] at (axis cs: 0.7,2.5) {$f(x)=x^{2}-3x^{1/3}+5/2$};
\node[style={fill=white}] at (axis cs: 0.7,2.0) {$\Delta_{\ast}=0.0481$};

\end{axis}

\end{tikzpicture}
\caption{Approximation of the function $f(x)$ by a polynomial $P_{\ast}(x)$ with $N=7$ terms.}
\label{F-AFPN7}
\end{figure}
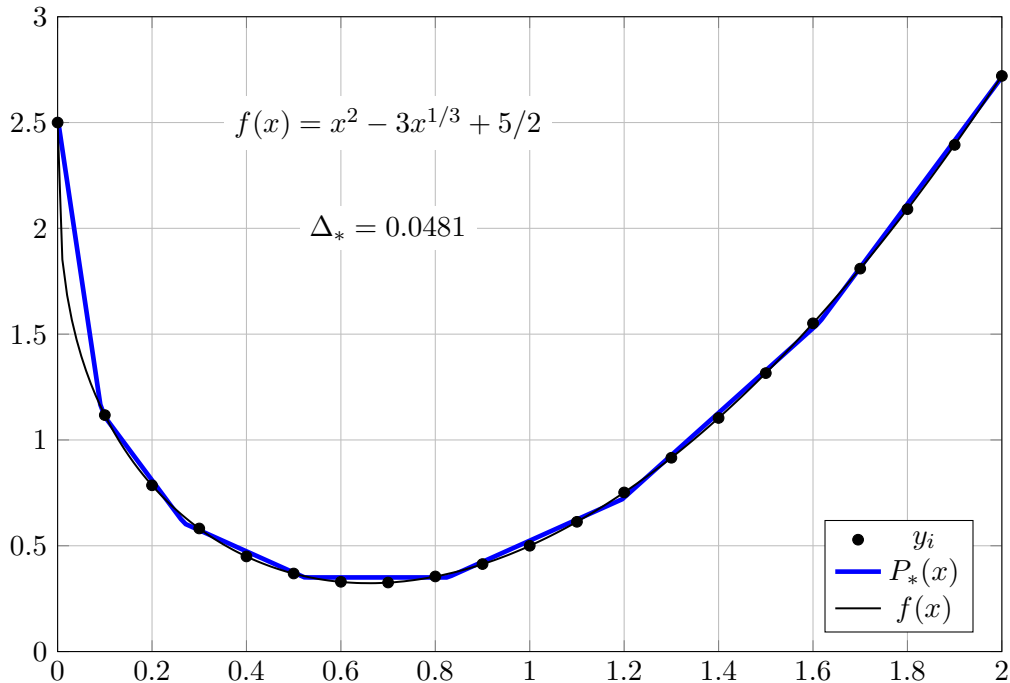

The obtained solution takes the form of the function 
\begin{multline*}
P_{\ast}(x)
=
\max(2.5240-15x,
1.4096-3x,
0.8736-x,
0.3503,
\\
-0.4760+x,
-1.6720+2x,
-3.2853+3x
),
\end{multline*}
which is demonstrated in Figure~\ref{F-AFPN7}.

We observe that, as expected, the approximation by polynomials with the number of terms $N$ increases from $5$ to $7$ results in a reduction in approximation error, which decreases the squared error $\Delta_{\ast}$ from $0.1360$ to $0.0481$.

\subsection{Discrete Approximation of Nonconvex Function}

Once again, we suppose that $M=21$ samples $(x_{i},y_{i})$ are given, which are obtained at the points $x_{i}=(i-1)/10$ for $i=1,\ldots,M$, from the function
\begin{equation*}
g(x)
=
3(x-1)^{2}\sin(x)+1/4,
\qquad
0\leq x\leq 2.
\end{equation*}

This function is nonconvex, which makes it appropriate to approximate the sample data by a tropical rational function defined as the ratio of polynomials:
\begin{equation*}
R(x)
=
P(x)/Q(x).
\end{equation*}

We take numerator and denominator polynomials with $N=4$ and $L=2$ terms to be
\begin{equation*}
P(x)
=
\theta_{1}x^{p_{1}}\oplus\cdots\oplus\theta_{4}x^{p_{4}},
\qquad
Q(x)
=
\sigma_{1}x^{q_{1}}\oplus\sigma_{2}x^{q_{2}}.
\end{equation*}

Since the vectors $\bm{p}=(p_{1},\ldots,p_{4})^{T}$ and $\bm{q}=(q_{1},q_{2})^{T}$ are not known, we obtain their components by random sampling from the discrete uniform distribution over $[-10,10]$. For each obtained pair of vectors, we employ the procedure of rational approximation described in the previous section.

Solving the approximation problem for 10,000 sample vector pairs results in the minimum squared error $\Delta_{\ast}=0.1395$. The optimal vectors of coefficients, together with corresponding vectors of degrees, are found in the following form:
\begin{gather*}
\bm{\theta}_{\ast}^{T}
=
\left(
\begin{array}{cccc}
11.4688 & 9.6386 & 5.3845 & -10.5453
\end{array}
\right),
\qquad
\bm{\sigma}_{\ast}^{T}
=
\left(
\begin{array}{cc}
11.1281 & 10.3861
\end{array}
\right),
\\
\bm{p}_{\ast}^{T}
=
\left(
\begin{array}{cccc}
-3 & -2 & 1 & 2
\end{array}
\right),
\qquad
\bm{q}_{\ast}^{T}
=
\left(
\begin{array}{cc}
-5 & -2
\end{array}
\right).
\end{gather*}

In the ordinary notation, the approximating function is written as
\begin{equation*}
R_{\ast}(x)
=
P_{\ast}(x)-Q_{\ast}(x),
\end{equation*}
where
\begin{align*}
P_{\ast}(x)
&=
\max(3.4753-3x, 2.7409-2x, -1.0110+x, -2.6014+2x),
\\
Q_{\ast}(x)
&=
\max(3.2525-5x,2.4211-2x).
\end{align*}

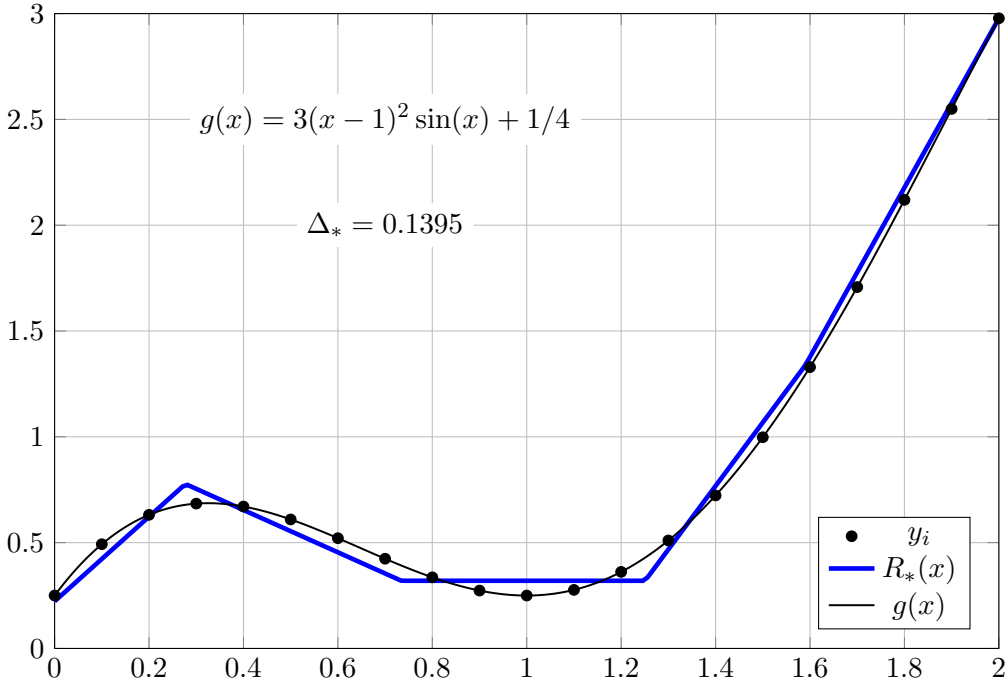
\begin{figure}[h]
\begin{tikzpicture}

\begin{axis}[legend pos=south east,
width=14cm,
height=10cm,
grid=major,
ymin=0,
ymax=3,
xmin=0,
xmax=2,
]

\addplot[
only marks,
]
coordinates {
(0,      0.2500)  
(0.1000, 0.4926)
(0.2000, 0.6314)
(0.3000, 0.6844)
(0.4000, 0.6706)
(0.5000, 0.6096)
(0.6000, 0.5210)
(0.7000, 0.4239)
(0.8000, 0.3361)
(0.9000, 0.2735)
(1.0000, 0.2500)
(1.1000, 0.2767)
(1.2000, 0.3618)
(1.3000, 0.5102)
(1.4000, 0.7230)
(1.5000, 0.9981)
(1.6000, 1.3295)
(1.7000, 1.7077)
(1.8000, 2.1198)
(1.9000, 2.5495)
(2.0000, 2.9779)
};

\addlegendentry{$y_{i}$}

\addplot[
blue,
samples=200,
line width=1.75pt,
domain=0.0:2.0,
y domain=0:10,
]{max(
3.4753+x*(-3),
2.7409+x*(-2),
-1.0110+x*1,
-2.6014+x*2
)
-max(
3.2525+x*(-5),
2.4211+x*(-2)
)};

\addlegendentry{$R_{\ast}(x)$}

\addplot[
samples=200,
black,
line width=0.75pt,
domain=0.0:2.0,
y domain=0:10,
]{3*(x-1)^2*sin(deg(x))+1/4};

\addlegendentry{$g(x)$}

\node[style={fill=white}] at (axis cs: 0.7,2.5) {$g(x)=3(x-1)^{2}\sin(x)+1/4$};
\node[style={fill=white}] at (axis cs: 0.7,2.0) {$\Delta_{\ast}=0.1395$};

\end{axis}

\end{tikzpicture}
\caption{Approximation of function $g(x)$ by a rational function $R_{\ast}(x)$, with $N=4$ and $L=2$.}
\label{F-AFRN4_2}
\end{figure}

A graphical illustration is given in Figure~\ref{F-AFRN4_2}.

To obtain a more accurate solution, we consider an approximating rational function that has numerator and denominator polynomials with $N=6$ and $L=4$ terms given by
\begin{equation*}
P(x)
=
\theta_{1}x^{p_{1}}\oplus\cdots\oplus\theta_{6}x^{p_{6}},
\qquad
Q(x)
=
\sigma_{1}x^{q_{1}}\oplus\cdots\oplus\sigma_{4}x^{q_{4}}.
\end{equation*}

After random sampling of 10,000 vector pairs, to define degrees of polynomials involved, and solving the corresponding rational approximation problems, we obtain a solution with squared error $\Delta_{\ast}=0.0701$ and the following vectors:
\begin{gather*}
\bm{\theta}_{\ast}^{T}
=
\left(
\begin{array}{cccccc}
6.9455 & 6.0860 & 4.9978 & 3.7461 & 0.7639 & -2.6361
\end{array}
\right),
\\
\bm{\sigma}_{\ast}^{T}
=
\left(
\begin{array}{cccc}
6.6880 & 6.2962 & 5.8009 & 2.4211
\end{array}
\right),
\\
\bm{p}_{\ast}^{T}
=
\left(
\begin{array}{cccccc}
-3 & -2 & 1 & 0 & 2 & 4 
\end{array}
\right),
\qquad
\bm{q}_{\ast}^{T}
=
\left(
\begin{array}{cccc}
-5 & -3 & -2 & 0
\end{array}
\right).
\end{gather*}

For the approximating function in the usual form $R_{\ast}(x)=P_{\ast}(x)-Q_{\ast}(x)$, we obtain
\begin{align*}
P_{\ast}(x)
&=
\max(6.9455-3x, 6.0860-2x, 4.9978-x, 3.7461, 0.7639+2x, -2.6361+4x),
\\
Q_{\ast}(x)
&=
\max(6.6880-5x, 6.2962-3x, 5.8009-2x, 2.4211).
\end{align*}

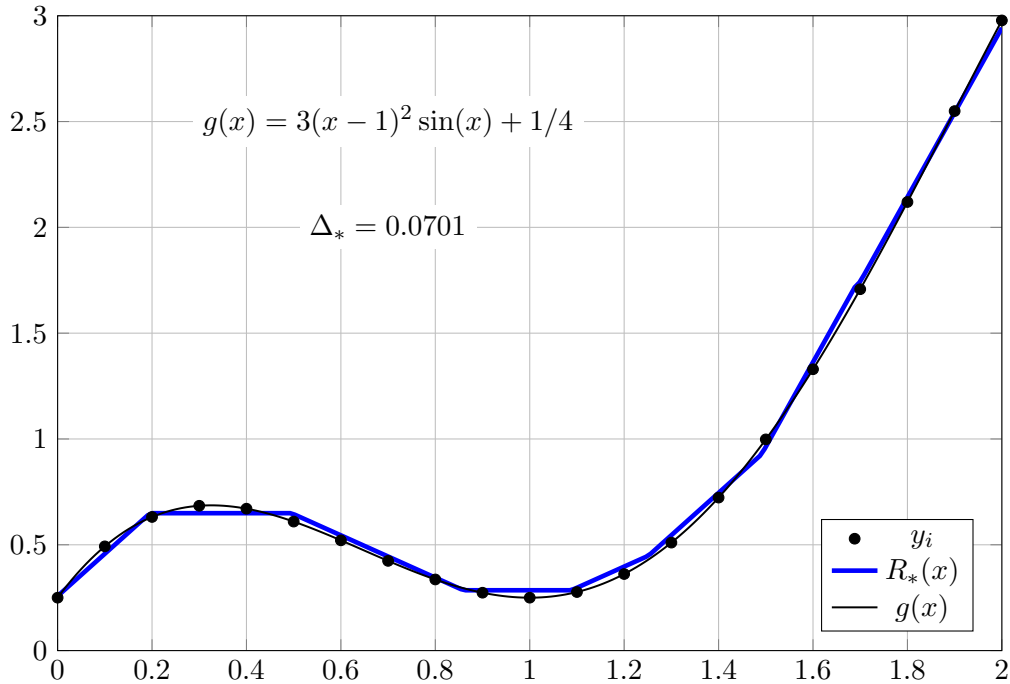
\begin{figure}[h]
\begin{tikzpicture}

\begin{axis}[legend pos=south east,
width=14cm,
height=10cm,
grid=major,
ymin=0,
ymax=3,
xmin=0,
xmax=2,
]

\addplot[
only marks,
]
coordinates {
(0,      0.2500)  
(0.1000, 0.4926)
(0.2000, 0.6314)
(0.3000, 0.6844)
(0.4000, 0.6706)
(0.5000, 0.6096)
(0.6000, 0.5210)
(0.7000, 0.4239)
(0.8000, 0.3361)
(0.9000, 0.2735)
(1.0000, 0.2500)
(1.1000, 0.2767)
(1.2000, 0.3618)
(1.3000, 0.5102)
(1.4000, 0.7230)
(1.5000, 0.9981)
(1.6000, 1.3295)
(1.7000, 1.7077)
(1.8000, 2.1198)
(1.9000, 2.5495)
(2.0000, 2.9779)
};

\addlegendentry{$y_{i}$}

\addplot[
blue,
samples=200,
line width=1.75pt,
domain=0.0:2.0,
y domain=0:10,
]{max(
6.9455+x*(-3),
6.0860+x*(-2),
4.9978+x*(-1),
3.7461+x*(0),
0.7639+x*(2),
-2.6361+x*(4)
)
-max(
6.6880+x*(-5),
6.2962+x*(-3),
5.8009+x*(-2),
2.4211+x*(0)
)};

\addlegendentry{$R_{\ast}(x)$}

\addplot[
samples=200,
black,
line width=0.75pt,
domain=0.0:2.0,
y domain=0:10,
]{3*(x-1)^2*sin(deg(x))+1/4};

\addlegendentry{$g(x)$}

\node[style={fill=white}] at (axis cs: 0.7,2.5) {$g(x)=3(x-1)^{2}\sin(x)+1/4$};
\node[style={fill=white}] at (axis cs: 0.7,2.0) {$\Delta_{\ast}=0.0701$};

\end{axis}

\end{tikzpicture}
\caption{Approximation of function $g(x)$ by a rational function $R_{\ast}(x)$, with $N=6$ and $L=4$.}
\label{F-AFRN6_4}
\end{figure}

The solution obtained is illustrated in Figure~\ref{F-AFRN6_4}.

\section{Conclusions}
\label{S-C}

In this paper, we introduced new discrete best approximation problems given in the tropical algebra setting. The problems were formulated to approximate an unknown function on a tropical (idempotent) semifield by tropical Puiseux polynomial and rational functions. We proposed an approach to solving both polynomial and rational approximation problems, by finding the best approximate solutions of the corresponding tropical linear vector equations. Specifically, the approximation by polynomials was reduced to a one-sided vector equation, for which we provided a direct analytical solution. The rational approximation led to a two-sided equation, which we solved using an iterative alternating algorithm. 

To illustrate the solution approach, we gave examples of solving approximation problems in terms of max-plus algebra. If defined in the context of max-plus algebra, Puiseux polynomials take the form of convex piecewise linear functions, whereas rational functions take the form of nonconvex piecewise linear functions. For the sake of simplicity, we used approximating polynomials with integer exponents that, in general, provide less flexibility of fitting data than polynomials with rational exponents. Nevertheless, the examples clearly demonstrated the applicability of the technique, at least to solve problems of discrete best Chebyshev approximation of real functions.

Moreover, we observed that the solution procedures proposed to solve the approximation problems did not change when the Puiseux polynomials were replaced by polynomials that may have had real exponents (if defined in the underlying semifield) and not just rational ones. This extended the approach, to allow generalized tropical Puiseux polynomial and rational functions to also be used as approximating functions.

We believe the main purpose of the paper is to concentrate on the key component of the new approximation approach and, therefore, we have left for future research many issues that may arise in implementation. Specifically, for approximation problems in the context of max-plus algebra, it is important to study how the approximate solution and its error depend on the number of monomials in the approximating polynomial (rational) functions and on the values of their exponents. The results of the study could help in the development of optimization methods to find exponents that minimize the approximation error, which presents another line of future research. Finally, implementation of the technique --- in terms of other tropical semifields, including max-algebra --- and the analysis of related results are also of considerable interest.

\section*{Acknowledgments}
The author is very grateful to the anonymous reviewers for their valuable comments and suggestions, which have been incorporated in the revised manuscript.

\bibliographystyle{abbrvurl}

\bibliography{Algebraic_solution_of_tropical_best_approximation_problems}

\end{document}